\documentclass[conference]{IEEEtran}
\IEEEoverridecommandlockouts
\usepackage{cite}
\usepackage{amsmath,amssymb,amsfonts}
\usepackage{algorithmic}
\usepackage{graphicx}
\usepackage{textcomp}
\usepackage{xcolor}
\usepackage{bm}
\usepackage{subcaption}
\usepackage{array}
\usepackage{booktabs}
\def\BibTeX{{\rm B\kern-.05em{\sc i\kern-.025em b}\kern-.08em
    T\kern-.1667em\lower.7ex\hbox{E}\kern-.125emX}}

\newcommand{\tabref}[1]  {table~\ref{#1}}
\newcommand{\secref}[1]  {section~\ref{#1}}

\begin{document}

\title{A Bayesian perspective on classical control
\thanks{MB is a JSPS International Research Fellow supported by a JSPS Grant-in-Aid for Scientific Research (No. 19F19809).}
}

\author{\IEEEauthorblockN{Manuel Baltieri}
\IEEEauthorblockA{
\textit{Laboratory for Neural Computation and Adaptation, RIKEN Centre for Brain Science}, Wako, Saitama, Japan \\
manuel.baltieri@riken.jp}
}

\maketitle

\begin{abstract}
The connections between optimal control and Bayesian inference have long been recognised, with the field of stochastic (optimal) control combining these frameworks for the solution of partially observable control problems. In particular, for the linear case with quadratic functions and Gaussian noise, stochastic control has shown remarkable results in different fields, including robotics, reinforcement learning and neuroscience, especially thanks to the established \emph{duality} of estimation and control processes. Following this idea we recently introduced a formulation of PID control, one of the most popular methods from classical control, based on active inference, a theory with roots in variational Bayesian methods, and applications in the biological and neural sciences. In this work, we highlight the advantages of our previous formulation and introduce new and more general ways to tackle some existing problems in current controller design procedures. In particular, we consider 1) a gradient-based tuning rule for the parameters (or gains) of a PID controller, 2) an implementation of multiple degrees of freedom for independent responses to different types of signals (e.g., two-degree-of-freedom PID), and 3) a novel time-domain formalisation of the performance-robustness trade-off in terms of tunable constraints (i.e., priors in a Bayesian model) of a single cost functional, variational free energy.
\end{abstract}


\section{Introduction}
In the last few decades, the importance of probabilistic approaches to optimal control theory has been highlighted by different applications of Bayesian methods to problems of control. In his pioneering work, Bellman introduced Markov decision processes \cite{bellman1957markovian} as part of what is now known as stochastic optimal control \cite{aastrom1970introduction, anderson1990optimal, stengel1994optimal, todorov2006optimal, kappen2011optimal}. This formulation captured the intrinsic probabilistic nature of problems of optimal control and decision making, with state transitions, outcomes and actions/decisions that cannot always be easily described in purely deterministic terms. Bellman's approach extended on his own work on the dynamic programming method for (deterministic) optimal control, defining the Bellman equation, its uses and limitations, including the idea of the curse of dimensionality \cite{bellman1957dynamic}. Shortly after, Kalman introduced the notions of observability and controllability of a system \cite{kalman1960contributions}, with the former expressing the degree to which states can be estimated from noisy observations, and the latter representing the degree of control over a system when different manipulations are applied. Kalman also noticed that his filter was \emph{dual} to the linear quadratic regulator (LQR), a now well known method in (optimal) control theory that he also established \cite{kalman1960new}, showing how both solutions require solving Riccati equations, forward in time for filtering (error covariance matrix) and backward in time for control (Hessian of the cost-to-go function). In the following years, several results improved the treatment of stochastic optimal control problems, including for instance the separation principle \cite{wonham1968separation, stengel1994optimal}\footnote{Also known in econometrics as certainty equivalence property \cite{simon1956dynamic, theil1957note}, but see \cite{stengel1994optimal} for a possible distinction.} and its applications to the treatment of regulation in the presence of uncertainty, i.e., for (linear) partially observable control problems. Due to its analytical tractability and its combination of estimation and control algorithms, the linear quadratic framework (i.e., stochastic optimal control for linear state-space models with Gaussian white noise and quadratic cost functions) has since then become a standard approach in different fields, including not only control theory and engineering \cite{anderson1990optimal, stengel1994optimal}, but also robotics \cite{lewis2003robot} and neuroscience \cite{todorov2002optimal, todorov2005stochastic}.

In recent years, the results based on the notion of duality in the linear case have been extended to (some classes of) nonlinear systems \cite{mitter2003variational, kappen2005linear, todorov2008general}, highlighting further connections between control and estimation. Notably, these extended dualities often rely on more efficient variational approximations commonly used in problems of inference. For instance, relevant advances in stochastic optimal control and reinforcement learning have been driven by the use of methods commonly adopted to approximate intractable problems of Bayesian inference, e.g., variational Bayes. These methods have been shown to outperform standard dynamic programming and reinforcement learning algorithms for the control of different classes of problems \cite{attias2003planning, mitter2003variational, todorov2009efficient, kappen2011optimal}. Building on these ideas, a similar approach has been proposed and adopted in neuroscience in an attempt to characterise brain function and sensorimotor control under a unifying probabilistic framework: active inference. While a full treatment of this framework is beyond the scope of the present work (for some technical reviews see, for instance, \cite{buckley2017free, da2020active}), we highlight how active inference combines methods from machine learning (variational Bayes), control theory (stochastic control) and statistical inference (hierarchical and empirical Bayes) to form a theory that includes several existing results from different fields as special cases, from predictive coding, to the infomax principle, to statistical models of learning, to risk-sensitive and KL-control \cite{Friston2008a, Friston2008c, Friston2010nature, friston2017active, baltieri2019active, da2020active}.

Most of these results rely, at the moment, on the application of (approximate) Bayesian approaches to \emph{optimal} control, with almost no mention of \emph{classical} control methods. While classical methods can be seen as a special case of optimal control, the possible advantages specific to Bayesian formulations of classical algorithms such as Proportional-Integral-Derivative (PID), remain largely unexplored. In this work we look at classical controllers from the perspective of approximate Bayesian inference, discussing the implications of variational Bayesian methods for the future of, in particular, PID control \cite{aastrom2001future}. This perspective has previously been adopted in, for instance, \cite{baltieri2019pid} where a new gradient-based gain tuning rule was derived in closed-form for optimal regulation near the set-point/reference goal. In the next sections we present three cases in support of a new (Bayesian) framework to design and study classical control methods that ought to be seen as complementary to existing ones, e.g., optimal control and frequency-domain analysis. We will first, briefly, 1) recapitulate the previously derived results for gain tuning introducing new connections to path integral control and estimation in the presence of biases, then 2) present a more formal and in depth treatment of the connections between PID controllers with two degrees of freedom and Bayesian inference schemes, and finally 3) focus on the open challenge of framing different competing constraints of the performance-robustness trade-off in PID control, here defined in terms of priors and hyperpriors on a probabilistic generative model.

\section{Case 1: A Bayesian derivation of PID gains and their optimisation}
PID controllers are the most popular choice for SISO systems regulation in different areas of engineering \cite{aastrom2010feedback}. Their popularity is mainly due to their simplicity and relatively low number of tunable parameters. However, despite only including a few key parameters, or gains, their tuning (or optimisation) remains largely an open challenge \cite{aastrom2001future, ang2005pid}. Existing tuning methods are often limited to specific cases or applications, relying on (ad-hoc) analytical rules, simple heuristics, frequency domain analysis, optimisation (including the use of artificial neural networks) or a combination of the above (for a survey, see \cite{ang2005pid}), that hardly generalise across different classes of problems. Here we report a more general method than can be explicitly derived by taking a different, Bayesian perspective on control problems. 

Previous work relating classical control to optimal observers, and thus indirectly to Bayesian methods, showed that the integral component of PID controllers corresponds to a process of estimation of unknown (but linear/constant or step) perturbations, equivalent to a Kalman-Bucy filter with augmented state for the inference of unknown inputs (or biases) \cite{johnson1968optimal, johnson1975observers}. Using the same approach, this connection was then generalised to higher order polynomial disturbances, equivalent to controllers including further integration terms \cite{johnson1970further}, i.e., corresponding to PIID, PIIID, etc. controllers. In \cite{baltieri2019pid}, we derived a fully probabilistic version of PID control, highlighting in particular some of the relationships between integral control and an emerging framework in computational and cognitive neuroscience, active inference \cite{Friston2010nature, friston2017active}. Using active inference, we thus defined a more explicit generative model to describe an underlying stochastic process producing PID control as a gradient descent of a cost functional, variational free energy. While these two approaches, \cite{johnson1970further} and \cite{baltieri2019pid}, share a number of features, they also present some core technical differences. Our proposal in fact includes:
\begin{itemize}
  \item a more direct interpretation of the control matrix $R$, commonly used as a weight for the cost of control in the value (or cost-to-go) function \cite{anderson1990optimal, stengel1994optimal},
  \item a gradient-based algorithm to optimise $R$, and
  \item a generalisation to (some classes of) nonlinear problems.
\end{itemize}
As shown in \cite{johnson1968optimal, johnson1970further}, the control matrix $R$ is particularly relevant for the computation of the gains of PID controllers, here treated as part of the feedback matrix of a linear quadratic controller. In active inference, such gains correspond to specific hyperparametrisations of the linear state space (generative) model used to approximate the dynamics of the system to control, i.e., the (expected) precision, or inverse covariance, of the observation noise \cite{baltieri2019pid}. This result is closely related to Kalman's duality of inference and control \cite{kalman1960general, kalman1960contributions, todorov2008general}, highlighting the mathematical correspondence between the processes of stochastic estimation and deterministic control, At the same time, the active inference formulation extends this duality beyond simply noting mathematical similarities, in order to include an account of the \emph{dual role} of action in the context of exploration/exploitation problems \cite{bar1974dual, baltieri2020kalman}. Furthermore, given the role of $R$ as (expected) precision in the generative model, the gain parameters of PID controllers can be optimised via a gradient descent on the same cost functional, i.e., variational free energy \cite{baltieri2019pid}, following a second order scheme introduced in \cite{Friston2010genfilt} that, under some assumptions, holds also for some classes of nonlinear problems \cite{Friston2008c, Friston2010genfilt}. 

\section{Case 2: PID control with 2DOF in active inference}
\label{sec:2DOFPID}
In many applications of PID control, it is often desirable to build regulators that can respond to external disturbances while avoiding large fluctuations (e.g., overshooting) due to changes of the target of the regulation process. In standard PID control, these requirements are shown to be conflicting \cite{aastrom1995pid, araki2003two} thus leading, in the most general case, to a multi-objective optimisation problem whose solutions lie on a Pareto front defined by
\begin{itemize}
  \item changes in the control target (i.e., set-point response), and 
  \item changes in the amplitude of a step disturbance (i.e., disturbance response).
\end{itemize}
To overcome the limitations induced by this trade-off, previous work (see \cite{araki2003two, aastrom2001future} and references therein) introduced the idea of controllers with two degrees of freedom, or 2DOF, PID. Multiple degrees of freedom obtained by augmenting controllers with multiple internal loops of PI or PID control (see also equivalent examples such as PI-PD control \cite{aastrom2001future}), then ensure that different constraints can be treated independently, using parameters from different sub-loops to encode separate desired behaviours \cite{araki2003two, aastrom2001future}.

Our probabilistic derivation of PID control via a variational approximation of Bayesian inference showcases a clear and direct interpretation of the presence of two degrees of freedom, here derived using rather general arguments. Unlike previous proposals, one need not augment a controller with an extra feedforward component that can separate the effects of a compensator for disturbances or set-point changes \cite{araki2003two}. In active inference, the existence of two degrees of freedom is a simple consequence of the probabilistic (Bayesian) description of the generative model used to derive a controller \cite{baltieri2019pid}. This becomes more obvious after looking at the variational free energy (see equation (13) in \cite{baltieri2019pid}) here reported as\footnote{Some terms in the free energy functional are hereby dropped for clarity. For the treatment of an extra set of terms important for the optimisation of PID gains, see \cite{baltieri2019pid}. For a more complete discussion of other terms which are constant during the optimisation phase, see \cite{Friston2008a, Friston2010genfilt, buckley2017free}.}
\begin{align}
    F \approx \frac{1}{2} \bigg[ \mu_{\pi_{\tilde{z}}} \Big(\tilde{y} - g(\tilde{\mu}_x, \tilde{\mu}_v)\Big)^2 + \mu_{\pi_{\tilde{w}}} \Big (\tilde{\mu}'_x - f(\tilde{\mu}_x, \tilde{\mu}_v) \Big)^2 \bigg]
  \label{eq:freeEnergyLaplace}
\end{align}
where $y, \mu_x, \mu_v$ are observations (or measurements), and expected hidden states (the estimate of the state of the system to regulate) and inputs (the set-point) respectively. Hyperparameters $\mu_{\pi_{z}}, \mu_{\pi_{w}}$ are the expected precisions on observation and system noise, and $f(), g()$ are state transition and observation functions. The tilde simply highlights a notation used to group different derivatives, or rather embeddings orders, of a variable, e.g., $\tilde{y} = \{y, y', y''\}$, see \cite{baltieri2019pid} for more details. The simplified (i.e., under Gaussian assumptions) variational free energy functional in \eqref{eq:freeEnergyLaplace} contains two sets of prediction errors, essentially instantiating two degrees of freedom for the controller. Notice that unlike equation (13) in \cite{baltieri2019pid}, here we explicitly replaced $\pi_{\tilde{z}}, \pi_{\tilde{w}}$ with $\mu_{\pi_{\tilde{z}}}, \mu_{\pi_{\tilde{w}}}$ from the beginning, to highlight the fact that hyperparameters $\mu_{\pi_{\tilde{z}}}, \mu_{\pi_{\tilde{w}}}$ are only \emph{estimates} of some ``true'' hyperparameters $\pi_{\tilde{z}}, \pi_{\tilde{w}}$. This follows from a full Bayesian treatment of the control problem, considering all variables to be \emph{random} variables \cite{robert2007bayesian} (cf. traditional point-estimates in frequentist frameworks for statistical learning). In our case, to simplify the mathematical treatment, we treat them as Gaussian random variables with means $\mu_{\pi_{\tilde{z}}}, \mu_{\pi_{\tilde{w}}}$ (and covariances to be discussed in the next section). Importantly, these expectations are updated on a slower time scale \cite{baltieri2019pid}, following schemes found in \cite{Friston2008a} and in particular \cite{Friston2010genfilt}, emphasising how parameters and hyperparameters of a generative model ought to be considered as fixed quantities over a certain (i.e., long) time scale.

The two sets of prediction errors, $\mu_{\pi_{\tilde{z}}} \Big(\tilde{y} - g(\tilde{\mu}_x, \tilde{\mu}_v)\Big)$ and $\mu_{\pi_{\tilde{w}}} \Big (\tilde{\mu}'_x - f(\tilde{\mu}_x, \tilde{\mu}_v) \Big)$ weighted by hyperparameters $\mu_{\pi_{\tilde{z}}}$ and $\mu_{\pi_{\tilde{w}}}$, represent likelihood and prior of a Bayesian update scheme formulated using generative models under Gaussian assumptions, the Laplace \cite{mackay2003information} and the variational Gaussian \cite{opper2009variational} approximations (to clarify their role see discussion in Chapter 3 of \cite{baltieri2019active}). The update equations minimising these prediction errors \cite{baltieri2019pid} (also called recognition dynamics \cite{kim2018recognition}), are similar to the update and prediction steps of standard algorithms from estimation theory such as Kalman-Bucy filters \cite{jazwinski1970stochastic}, and equivalent to feedback and feedforward loops in 2DOF PID controllers \cite{araki2003two}. In this set up, PID control with a single degree of freedom can be derived as the limit case for fully observable states, i.e., $\tilde{y} = \tilde{x}$ (cf. state feedback methods in \cite{aastrom2010feedback}). The independence of set-point and disturbance responses crucial for 2DOF PID controllers then corresponds, in this framework, to a generative model having system and measurement noise independent of one another, a standard assumption for linear state space models.

\section{Case 3: The performance-robustness trade-off for PID controllers under active inference}
\label{PerfRobTradeOff}
The presence of conflicting criteria for the design of PID controller is a well known issue in the control theory literature \cite{rivera1986internal}, as partially highlighted in the previous section. This conflict is often referred to as the performance-robustness trade-off \cite{aastrom2001future, aastrom2006advanced, garpinger2014performance}. Controllers are usually designed to optimise some given performance criteria while, at the same time, attempting to maintain a certain level of robustness in face of uncertainty and unexpected conditions during the regulation process. The performance of a controller is normally assessed using one or more of the following criteria \cite{aastrom2001future, aastrom2006advanced}:
\begin{itemize}
  \item \emph{load disturbance response}, or how a controller reacts to changes in external inputs, e.g., a step input,
  \item \emph{set-point response}, or how a controller responds to different set-points over time,
  \item \emph{measurement noise response}, or how noise on observations impacts the regulation process,
\end{itemize}
while robustness is mainly evaluated based on:
\begin{itemize}
  \item \emph{robustness to model uncertainty}, or how uncertainty on plant/environment dynamics affects the controller.
\end{itemize}
The goal of a general methodology for the design and tuning of PID controllers is to bring together these (and possibly more) criteria into a formal, unified and tractable framework that can be applied to a large class of compensation problems. An example in this direction is presented in \cite{aastrom1998design} (see also \cite{zhuang1993automatic, grimble1999algorithm, o2008optimal} for other partial attempts). This methodology is based on the maximisation of the integral gain (equivalent, near the reference point, to the minimisation of the integral of the error from the set-point \cite{aastrom1995pid}), subject to constraints derived from a frequency domain analysis related to the Nyquist stability criterion applied to the controlled system \cite{aastrom1998design}. Here, we propose our Bayesian formulation as an alternative (and in many cases complementary) framework for the design of PID controllers that leverages the straightforward interpretation of the performance-robustness trade-off for PID controllers in terms of uncertainty parameters (i.e., hyperparameters, precisions or inverse covariances) in the variational free energy \cite{baltieri2019pid}. To highlight its potential, we discuss the four standard criteria listed above as part of performance-robustness trade-off to address what can be gained using a Bayesian perspective.

\subsection{Load disturbance response}
\label{sec:Load}
A classic design principle for PID controllers is based on the response of a controller to perturbations that drive a process away from the target value \cite{aastrom1995pid}. Random, zero-mean disturbances are commonly modelled as white Gaussian variables, and the parameters of the controller are simply tuned to reject such noise. Integral control then guarantees an appropriate response to step disturbances, equivalent to \emph{non}-zero-mean noise (or to a bias term \cite{johnson1975observers}), by accumulating and compensating for the ensuing steady-state error \cite{johnson1968optimal, francis1976internal, aastrom1995pid, sontag2003adaptation}. The load disturbance response is usually expressed in terms of a minimisation of the Integral Absolute Error (IAE) between the state of the system to regulate and its target:
\begin{align}
    IAE = \int_o^{\infty} \left| e(t) \right| \; \text{d}t
\end{align}
or approximated by the Integral Error (IE) for non-oscillating or oscillating but well-damped systems \cite{aastrom1995pid}:
\begin{align}
    IE = \int_o^{\infty} e(t) \; \text{d}t
\end{align}
The IE criterion is especially relevant because it gives a straightforward intuition of the role of integral gain since, under a few simplifying assumptions (including a system's initial state close to the target value), the IE is equal to the inverse of $k_i$ as $t \rightarrow \infty$ \cite{aastrom1995pid}. This implies that for large (theoretically, infinite) integral gains, the IE is minimised. While useful for its straightforward interpretation of this free parameter, practical and physical limitations often restrict the maximisation of the integral gain.

Our formulation builds on previous work showing how the use of integral control is optimal for unknown step perturbations applied to a system \cite{johnson1968optimal, sontag2003adaptation}. In statistical terms, the presence of such disturbances can be formally seen as a \emph{bias} term in an estimation process \cite{friedland1969treatment}, showing how rejecting (step) perturbations is equivalent to estimating biases \cite{johnson1975observers}. In active inference we can extend this (exact) result for linear systems and disturbances to nonlinear cases (not limited to polynomial perturbations as in \cite{johnson1970further}), where a more general (but often only approximate or limited to special classes of nonlinearities) duality of estimation and control is obtained using variational and path integral formulations \cite{mitter2003variational, todorov2008general}, or via probability integral transforms in the form of hierarchical generative models \cite{Friston2008c}.

Furthermore, in our (Bayesian) formulation we gain a second and arguably deeper intuition on the role of the integral gain, which is now explicitly represented as one of the expected precisions (or inverse covariance) of observations $\tilde{y}$, i.e., $\mu_{\pi_{z}}$, see \cite{baltieri2019pid}. This prescribes a simple and alternative way of understanding why the maximisation of $k_i$ is usually a good heuristic for regulation problems where PID control is used \cite{aastrom2001future}: maximising $k_i$ is in fact equivalent to minimising uncertainty on measurements $y$, by maximising (minimising) the expected precision $\mu_{\pi_{z}}$ (variance $\mu_{\sigma_{z}}$) of the measurements of the system to regulate. At the same time, this can also explain some of the limitations of this heuristic, discussed in the frequency domain for instance in \cite{aastrom1998design}. The maximisation of $k_i$, without any constraints, corresponds to the minimisation of the expected measurement variance $\mu_{\sigma_{z}}$, such that $t \to \infty$,
\begin{align}
  \mu_{\sigma_{z}} \to 0.
\end{align}
In practice, however, one should always consider a certain level of intrinsic, i.e., aleatoric, uncertainty whose variance is fundamentally irreducible. Even an optimal controller can't overcome the limited sensitivity of a sensor (here represented by the ``real'' $\sigma_{z}$, as opposed to its estimate $\mu_{\sigma_{z}}$), bringing $\mu_{\sigma_{z}}$ down to 0 is thus not possible if $\sigma_{z} > 0$. In other proposals \cite{aastrom1998design}, the same aleatoric uncertainty $\sigma_{z}$ is effectively approximated with a measure that captures the levels of controllability of a system through the definition of appropriate sensitivity functions in the frequency domain.

Our Bayesian implementation also extends the intuition behind the integral gain as a precision of observations to the other two gains, $k_p$ and $k_d$. In our formulation, these gains become in fact estimated precisions of higher embedding orders of the observations, $y', y''$, often also called \emph{generalised coordinates of motion} \cite{Friston2008a, Friston2010genfilt}. These embedding orders essentially represent a Taylor expansion (in time) of continuous random variables defined according to a Stratonovich (rather than Ito) interpretation, equivalent to non-Markovian (semi-Markovian, of Markovian of order \emph{n}) stochastic processes \cite{jazwinski1970stochastic, Friston2008c, baltieri2019pid}. In practice, for measurements taken at a high enough frequency, and with controllers having a short enough intrinsic time scale to regulate such high frequency measurements (i.e., a time scale approaching the underlying continuous models of the systems to regulate), observation noise should be treated as coloured, rather than white as in standard delta-autocorrelated noise. Under these assumptions, the implementation of PID control and its extensions (e.g., multiple I and D terms) becomes simply a linear approximation of a measured non-Markovian trajectory. Perhaps in a more intuitive way, we can see expected precisions $\mu_{\pi_{\tilde{z}}} = \{\mu_{\pi_{z}}, \mu_{\pi_{z'}}, \mu_{\pi_{z''}}\}$ as simultaneously 1) representing the precision of a trajectory in the state-space (rather than the precision on a point) and 2) regulating the convergence rate of measurements to a set-\emph{trajectory} (rather than point), specifying how quickly a controller ought to respond to a sudden change in a set of observations and their higher orders of motion.

\subsection{Set-point response}
Following the load disturbance rejection property, a second performance criterion used for the design of PID controllers is their set-point response, i.e., how controllers respond to variations in the set value used as a target to regulate a system. Naively, this could be seen as closely related to load disturbances: rather than changes in the measurement, we now have changes in the target value, both of them used to define some error term, $e$. In practice however, it is desirable to decouple these two problems, creating a controller with different sensitivities to load disturbances or set-point updates whenever necessary \cite{aastrom1995pid}. This requires a controller with two degrees of freedom, as discussed in more detail in \secref{sec:2DOFPID}, which is an inherent feature of the active inference formulation where expectations of hidden states $\tilde{\mu}_x$ are updated using a (Bayesian) scheme that balances (via a set of independent expected precisions, or weights) prediction errors on both
\begin{itemize}
  \item observations, $\mu_{\pi_{\tilde{z}}} \Big(\tilde{y} - g(\tilde{\mu}_x, \tilde{\mu}_v)\Big)$, where load disturbances can appear as part of the measurements $\tilde{y}$, and
  \item system dynamics, $\mu_{\pi_{\tilde{w}}} \Big (\tilde{\mu}'_x - f(\tilde{\mu}_x, \tilde{\mu}_v) \Big)$, where set-trajectories can be updated as inputs/priors $\tilde{\mu}_v$.
\end{itemize}
Mirroring the role of $\mu_{\pi_{\tilde{z}}}$ for load disturbances, expected precisions $\mu_{\pi_{\tilde{w}}}$ on dynamic prediction errors effectively implement a response mechanism to  set-trajectories updates, with high expected precisions implying a fast response, and low precisions entailing a slow one. Equivalently, from a probabilistic perspective this can be explained with the idea that the former describes a model with low uncertainty on dynamics (high precision = low covariance) meaning that any variation from such dynamics should quickly be dealt with; the latter encodes, on the other hand, the fact that high expected covariance allows for changes in set-trajectories, i.e., sudden updates are not surprising, therefore changes can be slow (and in the limit for very large covariances, almost absent).

\subsection{Measurement noise response}
A third common requirement for PID controllers is related to their performance in face of noisy or uncertain measurements. These may be due, for instance, to physical constraints/sensitivities of the sensors. In the literature, high frequency measurement noise \cite{aastrom2010feedback} is usually tackled via a careful and ad-hoc controller design, including for example pre-filtering of the observed data \cite{aastrom2006advanced}. In the Bayesian formulation of PID controllers that we introduced, we have a direct measure of (the best estimate of) the measurement noise: the expected precision or inverse covariance $\mu_{\pi_{\tilde{z}}}$ of the random variable $z$ and its higher orders of motion. Measurement noise is thus related to the same set of hyperparameters used to explain load disturbances rejection which, on the other hand, can be seen as low frequency noise. This shows another trade-off between design criteria, in this case related to the high frequency properties of measurement noise and the (usually) low frequency of external disturbances.

The previously identified maximisation of expected precision $\mu_{\pi_{z}}$ (integral gain $k_i$) implies an increased cutoff frequency of the low-pass filter implemented by linear generative models of the kind we introduced to approximate PID control \cite{andrews2006optimal}. This suggests that, while low-frequency disturbances can be suppressed more quickly (even if at the cost of possibly overshooting), this comes at the expenses of a ``hypersensitivity'' to high-frequency noise, i.e., not rejecting as much noise as otherwise possible with slower load disturbance responses (as shown in a simple model, for instance, in \cite{andrews2006optimal}). In our framework, this can easily be noticed by looking at the role played by expected precision $\mu_{\pi_{\tilde{z}}}$ in the time domain, encoding expected variability in observed data without a clear distinction between rare perturbations and persistent noise.

At the same time, however, the active inference formulation can be used to treat this problem in a more principled way, introducing informative priors on expected precisions $\mu_{\pi_{\tilde{z}}}$, i.e., hyperpriors $\eta_{\pi_{\tilde{z}}}$\footnote{To maintain a notation similar to the one used in \cite{Friston2008c, baltieri2019pid}.} (or more complicated functions $h(\eta_{\pi_{\tilde{z}}})$), see \cite{Friston2008c} for a formal treatment. The variational free energy functional then includes another set of prediction errors, (cf. \eqref{eq:freeEnergyLaplace}),
\begin{align}
    F \approx \frac{1}{2} \bigg[ & \mu_{\pi_{\tilde{z}}} \Big(\tilde{y} - g(\tilde{\mu}_x, \tilde{\mu}_v)\Big)^2 + \mu_{\pi_{\tilde{w}}} \Big (\tilde{\mu}'_x - f(\tilde{\mu}_x, \tilde{\mu}_v) \Big)^2 \nonumber \\
    + & \mu_{p_{\tilde{z}}} \Big (\mu_{\pi_{\tilde{z}}} - h(\eta_{\pi_{\tilde{z}}}) \Big)^2
    \bigg]
  \label{eq:freeEnergyLaplaceExtended1}
\end{align}
with $\mu_{p_{\tilde{z}}} \Big (\mu_{\pi_{\tilde{z}}} - h(\eta_{\pi_{\tilde{z}}}) \Big)$ playing the role of L2 (or Tikhonov) regularisation terms for the ensuing recognition dynamics derived as a gradient descent on \eqref{eq:freeEnergyLaplaceExtended1}. Using these prediction errors one can effectively encode, for instance, constraints that reject strong high frequency noise by specifically targeting frequent large instantaneous fluctuations of expected precision $\mu_{\pi_{\tilde{z}}}$, penalising them with a Gaussian hyperprior (an L2 regularisation term that affects ``outliers'') centred at $\eta_{\pi_{\tilde{z}}}$. While such hyperprior would certainly then also influence the response to sparse step changes, expected precisions $\mu_{\pi_{\tilde{z}}}$ could be updated by slowly shifting hyperpriors $\eta_{\pi_{\tilde{z}}}$ to reflect biases in measurements $\tilde{y}$ that persist over a long period of time. Importantly, while the cost functional presents in this case some new terms, the underlying minimisation scheme remains the same: the recognition dynamics will simply include extra regularisation terms while still following a gradient descent on the (augmented) variational free energy. 

In the same way expected precisions $\mu_{\pi_{\tilde{z}}}$ regulate the response to changes in the observations due to load disturbances, expected precisions on higher order stochastic properties (e.g., expected precisions on expected precisions, $\mu_{p_{\tilde{z}}}$) can then be seen as regulating how a controller adapts to varying levels of measurement noise covariance given some (informative) priors $h(\eta_{\pi_{\tilde{z}}})$. For example, in cases where the variance of measurement noise changes over time, e.g., due to the natural degradation of sensors, our formulation can include mechanisms that take into account existing prior knowledge and that can be used by a controller to dynamically adapt to new levels of noise. More in general, in the presence of \emph{stochastic volatility} (i.e., models where the covariances of different random variables are themselves random variables \cite{robert2007bayesian}), one can easily encode prior knowledge of higher order properties of random variables by including extra hierarchical layers on the generative model we introduced for PID control.

\subsection{Robustness to model uncertainty}
PID controllers are usually designed to withstand some level of model uncertainty, inherent in any system we observe, interact with and try to regulate. In control theory, this problem affects compensators attempting to regulate a system while having access only to a limited amount of information regarding the dynamics of the system itself. PID controllers are especially popular as a ``model free'' strategy, or rather, for the small number of tunable parameters that are necessary to afford robust, although often suboptimal, control \cite{rivera1986internal}. In control problems, this robustness is sometimes captured by sensitivity functions \cite{aastrom1998design,aastrom2006advanced}, providing a proxy for, among other things, the sensitivity of a feedback system to variations in models of process dynamics. In our derivation of PID control as a process of Bayesian (active) inference, the uncertainty of the dynamics is represented by the expected precisions of system dynamics, $\mu_{\pi_{\tilde{w}}}$, in the linear generative model defined in \cite{baltieri2019pid}. For instance, low expected precisions $\mu_{\pi_{\tilde{w}}}$, expressing high uncertainty/covariance, encode the (Bayesian) belief that large fluctuations in the dynamics can be expected, while high expected precisions express the fact that dynamics should show only small fluctuations. Moreover, using our formulation we can describe the behaviour of a PID controller such that under controllability assumptions \cite{kalman1960contributions, stengel1994optimal}, it effectively ``imposes'' its own (linear) dynamics/priors on a system through larger weighted prediction errors $\mu_{\pi_{\tilde{w}}} \Big (\tilde{\mu}'_x - f(\tilde{\mu}_x, \tilde{\mu}_v) \Big)$, by forcing it into an attractor encoded in the set-trajectory represented by the controller's priors. The state of affairs of the world is only partially relevant to a PID controller, since as long as conditions of reachability and controllability \cite{stengel1994optimal} are met, all it does is try to drive a (controllable) system towards the desired equilibrium encoded by its priors on a set-trajectory.

As in the case of measurement noise, our formulation allows for the construction of an extra layer of hyperpriors to handle model uncertainty: in the active inference formulation we can in fact include priors on expected precisions $\mu_{\pi_{\tilde{w}}}$ to represent existing information on the expected/desired dynamics of a system to regulate
\begin{align}
    F \approx \frac{1}{2} \bigg[ & \mu_{\pi_{\tilde{z}}} \Big(\tilde{y} - g(\tilde{\mu}_x, \tilde{\mu}_v)\Big)^2 + \mu_{\pi_{\tilde{w}}} \Big (\tilde{\mu}'_x - f(\tilde{\mu}_x, \tilde{\mu}_v) \Big)^2 \nonumber \\
    + & \mu_{p_{\tilde{w}}} \Big (\mu_{\pi_{\tilde{w}}} - k(\eta_{\pi_{\tilde{w}}}) \Big)^2
    \bigg]
  \label{eq:freeEnergyLaplaceExtended2}
\end{align}
For instance, it is not hard to imagine that, following standard hierarchical or empirical Bayes methods in statistical inference \cite{robert2007bayesian}, information on existing control problems could be used to define classes of systems whose shared statistical properties form generic priors $\eta_{\pi_{\tilde{w}}}$. These priors could then be used to initialise our model in a suitable part of the state space to ensure a desired level of robustness (and if a similar approach were to be adopted for $\eta_{\pi_{\tilde{z}}}$, to guarantee some desired performances). In such settings, expected precisions $\mu_{\pi_{\tilde{w}}}$ can still be optimised via a simple gradient descent, now L2-regularised with newly introduced priors entering the variational free energy equation in the form of weighted prediction errors, $\mu_{p_{\tilde{w}}} \Big (\mu_{\pi_{\tilde{w}}} - k(\eta_{\pi_{\tilde{w}}}) \Big)$. This approach, especially when employing empirical Bayes, is similar in spirit to the clever initialisation achieved in deep learning approaches via ``pre-training'', where introducing an unsupervised learning phase before supervised training showed substantial improvements in the performance and generalisation properties of neural networks \cite{harvey1996unicycling, lecun2015deep}.

\section{Discussion}
The duality of estimation and control has long been recognised and exploited in problems of regulation under constraints of partial observability, i.e., stochastic control \cite{kalman1960general, wonham1968separation, anderson1990optimal, stengel1994optimal, mitter2003variational, todorov2008general}. This property relies on the mathematical equivalence of some classes of estimation and regulation problems, formulated as Bayesian inference and optimal control respectively. The applications of this duality have led to a series of significant new results in different areas, such as reinforcement learning \cite{todorov2009efficient}, robotics \cite{theodorou2010generalized} and neuroscience \cite{todorov2002optimal} where methods of \emph{approximate} Bayesian inference are now often employed to improve existing solutions. In this work we built on some of these previous ideas, discussing possible applications of Bayesian inference theories and related approximations to methods from classical control. In particular, we focused on PID control and on our previous implementation of this method in terms of Bayesian \emph{active} inference \cite{baltieri2019pid}, proposing this as a general unifying framework for the design of PID controllers still largely missing to date \cite{aastrom2001future, ang2005pid}.

In \cite{baltieri2019pid} we recently introduced a gradient-based procedure for gain tuning, using an interpretation of these parameters as stochastic properties (i.e., expected precisions, or inverse variances) of the system to regulate. Here we expanded on this formulation by providing direct links to Kalman's duality \cite{kalman1960general, kalman1960contributions, todorov2008general} and Bayesian estimation in the presence of bias terms, i.e., unknown step inputs \cite{johnson1975observers}.

We then discussed standard problems such as the necessity of two degrees of freedom in order to afford independent responses to load and set-point changes \cite{araki2003two}. Using the probabilistic interpretation given in \cite{baltieri2019pid}, we then drew a comparison between a pragmatic introduction of two degrees of freedom \cite{araki2003two}, represented by feedback and feedforward sub-loops in standard 2DOF PID control, and the more principled formulation of active inference, aligned with update and prediction equations of filtering algorithms (e.g., Kalman-filters \cite{kalman1960new}), and the use of prior and likelihood densities in recursive Bayesian update schemes \cite{jazwinski1970stochastic}. 

Crucially, we then proposed to frame one of the major open challenges for methods like PID control, the general performance-robustness trade-off due to the presence of conflicting design criteria  \cite{aastrom2001future, ang2005pid}, in terms of variational free energy minimisation \cite{kappen2005linear, Friston2008a, Friston2010genfilt, theodorou2012relative, baltieri2019pid}\footnote{The following equation combines \eqref{eq:freeEnergyLaplaceExtended1} and \eqref{eq:freeEnergyLaplaceExtended2}.}
\begin{align}
    F \approx & \frac{1}{2} \bigg[ \mu_{\pi_{\tilde{z}}} \Big(\tilde{y} - g(\tilde{\mu}_x, \tilde{\mu}_v)\Big)^2 + \mu_{\pi_{\tilde{w}}} \Big (\tilde{\mu}'_x - f(\tilde{\mu}_x, \tilde{\mu}_v) \Big)^2 \nonumber \\
    + & \mu_{p_{\tilde{z}}} \Big (\mu_{\pi_{\tilde{z}}} - h(\eta_{\pi_{\tilde{z}}}) \Big)^2
    + \mu_{p_{\tilde{w}}} \Big (\mu_{\pi_{\tilde{w}}} - k(\eta_{\pi_{\tilde{w}}}) \Big)^2
    \bigg]
  \label{eq:freeEnergyLaplaceExtended3}
\end{align}
In this formulation, simple constraints (load disturbance response and set-point change response) can easily be mapped to first order weighting parameters on the mean estimates of the state of the system to regulate. More complex ones, on the other hand, (measurement noise response and robustness to model uncertainty) can be introduced in terms of stochastic volatility \cite{robert2007bayesian}, i.e., by treating second moments (expected precisions, or hyperparameters) as random variable having appropriate hyperpriors encoded in the generative model. This mapping provides an immediate understanding of different desired statistical properties of the system to govern (see \tabref{tab:generalFramework}), now summarised in Table \ref{tab:generalFramework}.

\begin{table}[h]
\centering
\caption[Active inference as a general framework for PID controllers.]{{Active inference as a general framework for PID controllers design (adapted from \cite{baltieri2019pid} and here extended).}}
\begin{tabular}{m{1.4cm} m{1.4cm} m{5cm} }
	\toprule
	\textbf{Criterion} & \textbf{Mapped to} & \textbf{Interpretation in Active Inference} \\ \midrule
	Load disturbance response & \centering $\mu_{\pi_{\tilde{z}}}$ & Expected inverse covariance of the observations (i.e., precision), with low covariance implying a fast response, and vice versa \\ \midrule
	Set-point change response & \centering $\mu_{\pi_{\tilde{w}}}$ & Two degrees of freedom derived from the presence of two sets of prediction errors, sensory and dynamics, mapping to likelihood and priors of a Bayesian inference process \\ \midrule
	Measurement noise response & \centering (priors on) $\mu_{\pi_{\tilde{z}}}$ & Direct mapping of measurement noise to inverse covariance of the observations (i.e., precision), with hyperpriors (priors on expected precisions) introduced to differentiate high frequency noise from low frequency disturbances
	\\ \midrule
	Robustness to model uncertainty & \centering (priors on) $\mu_{\pi_{\tilde{w}}}$ & Direct mapping of model uncertainty to expected covariances of system fluctuations, representing the hidden dynamics of the system to control, with hyperpriors that can describe initial knowledge of, for example, a class of similar regulation problems to facilitate the optimisation of states/parameters (similar to the role of unsupervised ``pre-training'' in deep learning \cite{lecun2015deep}) \\
	\bottomrule
\end{tabular}
\label{tab:generalFramework}
\end{table}

\section{Conclusion and future work}
In an influential paper, {\AA}str{\"o}m and H{\"a}gglund asked whether PID control can play a role in the future of control theory and engineering \cite{aastrom2001future}. Despite being the most used controller in industry, the emergence of more specialised and better performing methods over the years, such as model predictive control, cast doubts on its long term applications and uses. {\AA}str{\"o}m and H{\"a}gglund however argued that due to its combined effectiveness and simplicity, PID is likely to remain relevant for the foreseeable future, perhaps in conjunction with other methods. At the same time, they highlighted a series of existing problems and open challenges faced by PID, including a relatively limited number of theoretical results in areas such as gain tuning and general (PID) controller design. In this work we built on our previous formulation of PID control in terms of active inference, a modern theory combining stochastic control and probabilistic Bayesian inference under the umbrella of variational free energy minimisation, to propose new applications of Bayesian methods to PID controllers in order to establish a more general design framework. After introducing a new practical implementation of optimal gain tuning in \cite{baltieri2019pid}, here we extended our proposal highlighting the connections between different design principles for PID, from the importance of multiple degrees of freedom to optimal tuning with conflicting performance-robustness criteria. This framework gives an interpretation of a series of different constraints as first and second order properties of a generative model that generates a PID controller as a gradient-based minimisation of a single cost functional, variational free energy. In the future, we will focus on simulations testing the current proposal using standard control benchmarks and following a vast literature on Bayesian models (see \cite{robert2007bayesian, Friston2008c} and references therein). We will then also draw more direct connections to modern machine and reinforcement learning, combining the present work with methods from \cite{tschantz2019scaling}, where preliminary results based on these and other ideas are utilised in the field of deep reinforcement learning with large neural networks performing amortised inference.

\footnotesize
\bibliographystyle{IEEEtran}
\bibliography{IEEEabrv,AllEntries} 

\end{document}